\documentclass{amsart}

\usepackage{amsmath,amssymb,amsthm,latexsym}

\newtheorem{theorem}{Theorem}
\newcommand{\bt}{\begin{theorem}}
\newcommand{\et}{\end{theorem}}
\newtheorem{lemma}{Lemma}
\newcommand{\bl}{\begin{lemma}}
\newcommand{\el}{\end{lemma}}
\newtheorem{corollary}{Corollary}
\newcommand{\bc}{\begin{corollary}}
\newcommand{\ec}{\end{corollary}}
\newcommand{\bconj}{\begin{conjecture}}
\newcommand{\econj}{\end{conjecture}}
\newtheorem{problem}{Problem}
\newcommand{\bprob}{\begin{problem}}
\newcommand{\eprob}{\end{problem}}
\newcommand{\beq}{\begin{equation}}
\newcommand{\eeq}{\end{equation}}
\newcommand{\benum}{\begin{enumerate}}
\newcommand{\eenum}{\end{enumerate}}

\newtheorem*{problemNN}{Problem}
\newcommand{\bprobNN}{\begin{problemNN}}
\newcommand{\eprobNN}{\end{problemNN}}

\newcommand{\N}{\ensuremath{ \mathbf N }}

\newcommand{\R}{\ensuremath{\mathbf R}}

\newcommand{\C}{\ensuremath{\mathbf C}}

\newcommand{\mcb}{\ensuremath{ \mathcal B}}

\newcommand{\mcu}{\ensuremath{ \mathcal U}}
\newcommand{\mcv}{\ensuremath{ \mathcal V}}

\newcommand{\mcx}{\ensuremath{ \mathcal X}}

\newcommand{\mba}{\ensuremath{ \mathbf a}}
\newcommand{\mbb}{\ensuremath{ \mathbf b}}
\newcommand{\mbc}{\ensuremath{ \mathbf c}}

\newcommand{\mbw}{\ensuremath{ \mathbf w}}
\newcommand{\mbx}{\ensuremath{ \mathbf x}}
\newcommand{\mby}{\ensuremath{ \mathbf y}}

\newcommand{\bmat}{\left(\begin{matrix}}
\newcommand{\emat}{\end{matrix}\right)}
\newcommand{\bsmallmat}{\left(\begin{smallmatrix}}
\newcommand{\esmallmat}{\end{smallmatrix}\right)}

\title{$B_h$-sets of real and complex numbers} 
\author{Melvyn B.  Nathanson}
\address{Department of Mathematics\\Lehman College (CUNY)\\Bronx, NY 10468}
\email{melvyn.nathanson@lehman.cuny.edu}

\date{\today}

\subjclass[2000]{11B05, 11B13, 11B30,11B34, 11B75, 54E52}
\keywords{Sidon set, $B_h$-set, open dense subsets, Baire's theorem,   
additive number theory, combinatorial number theory}
\thanks{Supported in part by  PSC-CUNY Research Award Program grant 66197-00 54.}

\begin{document}

\begin{abstract}
Let $K = \R$ or \C. 
An $n$-element subset $A$ of $K$ is a $B_h$-set if every $w \in K$ 
has at most one representation 
as the sum of $h$ not necessarily distinct elements of $A$.  
 Associated to the $B_h$ set $A = \{a_1,\ldots, a_n\}$ are the 
$B_h$-vectors $\mba = (a_1,\ldots, a_n)$ in $K^n$. 
This paper proves that ``almost all'' $n$-element  subsets of $K$ 
are $B_h$-sets in the sense that the set of all $B_h$-vectors is a dense open 
 subset of $K^n$.  
\end{abstract}

\maketitle

\section{Sumsets and $B_h$-sets} 
Let $A$ be a nonempty subset of an additive abelian group or semigroup $G$.    
For every positive integer $h$, the $h$-fold \emph{sumset}  of $A$  is the set 
of all sums of $h$ not necessarily distinct elements of $A$:
\[
hA = \underbrace{A+\cdots + A}_{\text{$h$ summands}}
= \left\{ a'_1+\cdots + a'_h: a'_i \in A \text{ for all } i \in [1,h] \right\}. 
\]
For integers $h \geq 2$ and $n \geq 2$, let  
\[
\mcx_{h,n} = \left\{ \mbx = (x_1,\ldots, x_n) \in \N_0^n: \sum_{i=1}^n x_i = h \right\}. 
\] 
For all $\mbx = (x_1,\ldots, x_n) \in \mcx_{h,n}$ and $\mby = (y_1,\ldots, y_n) \in \mcx_{h,n}$, we have  
\beq      \label{Sidon:xy-sum}
\| \mbx - \mby\|_{\infty}   = \max\{| x_i - y_i | : i=1,\ldots, n\}  \leq h.
\eeq

Let $\mba = (a_1,\ldots, a_n) \in G^n$ be a vector with distinct coordinates.  
Associated to the vector \mba\ is the subset $A = \{a_1,\ldots, a_n\}$ 
of $G$ of cardinality $n$.   
For all $\mbx \in \mcx_{h,n}$ and $\mba \in G^n$, we define the dot product 
\[
\mbx\cdot \mba = \sum_{i=1}^n x_ia_i.  
\]
Then $\mbx\cdot \mba  \in G$.  
The $h$-fold sumset of  $A $ can  be written in the form  
\[
hA = \left\{ \mbx \cdot \mba : \mbx  \in \mcx_{h,n} \right\}. 
\]

The subset  $A$ of $G$ is called   a \emph{$B_h$-set} 
if every element of $G$ has at most one representation (up to permutation of the summands) 
as the sum of $h$ not necessarily distinct elements of $A$.    
$B_2$-sets are also called \emph{Sidon sets}.  
The study of $B_h$-sets of integers is a classical topic in combinatorial additive number theory
(cf.~\cite{bose-chow62}--\cite{ruzs93a}).

This paper studies $B_h$-sets of real and complex numbers.  Let $K = \R$ or \C.  
Associated to every vector $\mba = (a_1,\ldots, a_n) \in K^n$ with distinct coordinates 
is the subset $A = \{a_1, \ldots, a_n\}$ of $K$ of cardinality $n$.  
The vector  \mba\ is called a  
\emph{$B_h$-vector} if its associated set $A$ is a $B_h$-set.  
Let $\mcb_h$ be the set of $B_h$-vectors.  
If $\mba = (a_1, a_2, \ldots, a_n) \in \mcb_h$, then 
$\sigma \mba = \left( a_{\sigma(1)}, a_{\sigma(2)},\ldots, a_{\sigma(n)}  \right)  \in \mcb_h$ 
 for every permutation $\sigma$ of $\{1,\ldots, n\}$. 

The object of this paper is to prove that ``almost all'' $n$-element  subsets of $K$ 
are $B_h$-sets in the sense that the set $\mcb_h$ is an open 
and dense subset of $K^n$.

\section{Open and dense subsets}

\bt                                                 \label{Bhvector:theorem:open}
Let $K = \R$ or \C. 
The  set $\mcb_h$  is an open subset of $K^n$.
\et

\begin{proof}
Let $\mba = (a_1,\ldots, a_n) \in K^n$ be a $B_h$-vector.  
If $\mbx,\mby \in \mcx_{h,n}$ with $\mbx \neq \mby$, 
then $\mbx\cdot \mba \neq \mby \cdot \mba$ and so 
\[
0 < \Delta = \min \left\{ \left\|  (\mbx\ - \mby) \cdot \mba  \right\|_{\infty} :
 (\mbx,\mby)  \in \mcx_{h,n}^2 \text{ and }  \mbx \neq \mby \right\}.
\]
We shall prove that, for all vectors $\mbb \in K^n$  with  
\[
0 < \| \mbb \|_{\infty} < \frac{\Delta}{h}  
\]
the vector 
\[
\mba + \mbb = \mbc 
\]
is a $B_h$-vector and so the set of $B_h$-vectors contains the open 
ball with center at \mba\and radius $\Delta/h$.   

Let  $\mbx,\mby \in \mcx_{h,n}$ with $\mbx \neq \mby$. 
If $\mbx\cdot \mbc = \mby \cdot \mbc$, then 
\[
\mbx\cdot (\mba + \mbb)  = \mby\cdot (\mba + \mbb) 
\]
and 
\[
(\mbx\ - \mby) \cdot \mba = (\mby - \mbx) \cdot \mbb  
\]
Applying inequality~\eqref{Sidon:xy-sum}, we obtain 
\begin{align*}
 \Delta & \leq \|(\mbx\ - \mby) \cdot \mba \|_{\infty} 
= \| (\mby - \mbx) \cdot \mbb  \|_{\infty}  \\
& \leq \| \mby - \mbx \|_{\infty}  \| \mbb  \|_{\infty}  \\ 
& < h \left(\frac{\Delta}{h}   \right)  = \Delta
\end{align*}
which is absurd.  Therefore, $\mbx \neq \mby$ 
implies $\mbx\cdot \mbc \neq\mby \cdot \mbc$ and so \mbc\ is a $B_h$-vector. 
This completes the proof. 
\end{proof}

\bl           \label{Bhvector:lemma:small}
For all $\delta > 0$ there is a $B_h$-vector \mbb\ such that $\|\mbb\|_{\infty} < \delta$.
\el

\begin{proof}
If $\mbw = (w_1,\ldots, w_n)$ is any $B_h$-vector in $\R^n$ or $\C^n$ 
with associated $B_h$-set $W =  \{w_1,\ldots, w_n\}$, 
then, for every $\lambda \neq 0$, 
the ``contraction'' $\lambda\ast W = \{\lambda w_i: i=1,\ldots, n\}$ is a $B_h$-set 
and the corresponding vector $\lambda \mbw = (\lambda w_1,\ldots, \lambda w_n)$ 
is also a $B_h$-vector.  Choosing  
$0 < \lambda < \delta/\|\mbw\|_{\infty}$ and $\mbb = \lambda\mbw$ gives 
\[
\|   \mbb\|_{\infty}  = \| \lambda \mbw\|_{\infty}  = |\lambda | \|\mbw \|_{\infty}  
< \left( \frac{\delta}{\|\mbw\|_{\infty} } \right) \|\mbw\|_{\infty} = \delta.
\]
This completes the proof.
\end{proof}

\bt                                      \label{Bhvector:theorem:dense}
Let $K = \R$ or \C.  
The  set $\mcb_h$ is a  dense subset of $K^n$.
\et

\begin{proof}
Let $\mba  \in K^n$ be a vector that is not a $B_h$-vector.  
We shall prove that, for every $\varepsilon > 0$, there is a $B_h$-vector 
$\mbc \in K^n$ such that $\|\mbc - \mba\|_{\infty} < \varepsilon$.

We partition the set of pairs of distinct vectors in $\mcx_{h,n}$ into two 
disjoint finite sets as follows:
\[
\mcu = \{ (\mbx,\mby) \in \mcx_{h,n}^2: \mbx \neq \mby \text{ and } 
\mbx\cdot \mba \neq \mby \cdot \mba \}
\]
and 
\[
\mcv= \{ (\mbx,\mby) \in \mcx_{h,n}^2: \mbx \neq \mby \text{ and } 
\mbx\cdot \mba = \mby \cdot \mba\}. 
\]
The set \mcv\ is nonempty because \mba\ is not a $B_h$-vector.    
The set \mcu\ is nonempty because $n \geq 2$ and $(\mbx,\mby) \in \mcu$ 
with $\mbx = (h,0,0,\ldots, 0)$ and  $\mby = (0,h,0,\ldots, 0)$.
Then 
\[
0 < \Delta = \min \left\{ \left\| (\mbx - \mby) \cdot \mba \right\|_{\infty} :
 (\mbx,\mby) \in \mcu \right\}.
\]

By  Lemma~\ref{Bhvector:lemma:small}, 
for all $\varepsilon > 0$ there is a $B_h$-vector \mbb\  in $K^n$ such that 
\[
\|\mbb\|_{\infty} <  \min \left( \varepsilon,\frac{\Delta}{h} \right).
\]
Applying inequality~\eqref{Sidon:xy-sum}, for all pairs $(\mbx,\mby) \in  \mcx_{h,n}^2$, 
we have 
\begin{align*}  
 \| ( \mbx  - \mby )\cdot \mbb     \|_{\infty}  
& \leq  \| \mbx  - \mby  \|_{\infty}  \|\mbb    \|_{\infty}  < h  \left(  \frac{\Delta}{h} \right) = \Delta.
\end{align*} 
The vector 
\[
 \mba + \mbb = \mbc 
\]
satisfies 
\[
\|\mbc - \mba\|_{\infty} = \|\mbb\|_{\infty} < \varepsilon.
\]

We shall prove that   \mbc\ is a $B_h$-vector.  
Equivalently, we shall prove that $\mbx\cdot \mbc \neq \mby\cdot \mbc$ 
for all pairs $(\mbx, \mby) \in \mcx_{h,n}^2$ with $\mbx \neq \mby$.   

If $(\mbx,\mby) \in \mcu$, then $\mbx\cdot \mba \neq \mby \cdot \mba$ and 
\begin{align*} 
\| \mbx \cdot \mbc - \mby \cdot \mbc\|_{\infty} 
& = \| (\mbx \cdot \mba - \mby \cdot \mba ) 
+ ( \mbx \cdot \mbb - \mby \cdot \mbb  )   \|_{\infty} \\
& \geq  \| (\mbx - \mby ) \cdot \mba   \|_{\infty} -
\| ( \mbx - \mby) \cdot \mbb    \|_{\infty} \\
& >  \Delta  -  \Delta  = 0. 
\end{align*} 

If $(\mbx,\mby) \in \mcv$, then $\mbx \cdot \mba = \mby \cdot \mba$.  
Because $\mbx \neq \mby$ and \mbb\ is a $B_h$-vector, we have 
$\mbx \cdot \mbb \neq  \mby \cdot \mbb$.  It follows that 
\begin{align*} 
\| \mbx \cdot \mbc - \mby \cdot \mbc\|_{\infty} 
& = \| (\mbx \cdot \mba - \mby \cdot \mba ) 
+ ( \mbx \cdot \mbb - \mby \cdot \mbb  )   \|_{\infty} \\
& = \| \mbx \cdot \mbb - \mby \cdot \mbb  \|_{\infty} \\
& >  0. 
\end{align*} 
This completes the proof.  
\end{proof}

Combining Theorems~\ref{Bhvector:theorem:open} 
and~\ref{Bhvector:theorem:dense} gives the following result.

\bt 
Let $K = \R$ or \C.  
The  set $\mcb_h$ is a  dense open subset of $K^n$.
\et

\bt
Let $\mcb_h$ be the set of $B_h$-vectors in $K^n$. 
The set  
\[
\mcb_{\infty} = \bigcap_{h=1}^{\infty} \mcb_h
\]
is a dense subset of $K^n$.
\et

\begin{proof}
This is simply an application of Baire's theorem in functional analysis.
\end{proof}


\section{$B_h[g]$-sets}

Let  $g$ and $h$ be positive integers.  
The subset  $A$ of an additive abelian semigroup $G$  is called a \emph{$B_h[g]$-set} 
if every element of $G$ has at most $g$ representations 
as the sum of $h$ not necessarily distinct elements of $A$.  
The $B_2[1]$-sets are the $B_2$-sets.  
If $A$ is a $B_h[g]$-set, then $A$ is also a $B_h[g']$-set for all $g' \geq g$.  
In particular, every $B_h$-set is a $B_h[g]$-set.  

Let $K = \R$ or \C.  The vector $\mba = (a_1,\ldots, a_n) \in K^n$ 
is a \emph{$B_h[g]$-vector} if the set $A = \{a_1,\ldots, a_n\}$ is an $n$-element 
$B_h[g]$-set in $K$.  Every $B_h$-vector is a $B_h[g]$-vector.  
Because the set of $B_h$-vectors is dense in  $K^n$, 
it follows that  the set of $B_h[g]$-vectors is also dense in $K^n$.

\bprobNN
Is  the set of $B_h[g]$-vectors also open in $K^n$?
\eprobNN

\def\cprime{$'$} \def\cprime{$'$} \def\cprime{$'$}
\providecommand{\bysame}{\leavevmode\hbox to3em{\hrulefill}\thinspace}
\providecommand{\MR}{\relax\ifhmode\unskip\space\fi MR }
\providecommand{\MRhref}[2]{%
  \href{http://www.ams.org/mathscinet-getitem?mr=#1}{#2}
}
\providecommand{\href}[2]{#2}

\end{document}